\documentclass[12pt,a4paper]{amsart}
\usepackage{amssymb}
\usepackage{mathrsfs}

\headsep=1cm \numberwithin{equation}{section}
\newtheorem{thm}{Theorem}[section]
\newtheorem{cor}[thm]{Corollary}
\newtheorem{lem}[thm]{Lemma}
\newtheorem{prop}[thm]{Proposition}
\theoremstyle{definition}

\theoremstyle{remark}

\numberwithin{equation}{section}

\newcommand{\inte}{\mathbb Z}
\newcommand\Supp{\operatorname{Supp}}
\newcommand\Ass{\operatorname{Ass}}

\newcommand\Spec{\operatorname{Spec}}
\newcommand\Rad{\operatorname{Rad}}

\newcommand\Ext{\operatorname{Ext}}

\newcommand{\fa}{\frak a }
\newcommand{\fb}{\frak b }
\newcommand{\p}{\frak p }
\newcommand{\q}{\frak q }

\begin{document} \title [Faltings' local-global principle  of local cohomology]
{A new proof of Faltings' local-global principle  for the finiteness of local cohomology modules}
\author{Davood Asadollahi and Reza Naghipour$^*$}
\address{Department of Mathematics, University of Tabriz, Tabriz, Iran;
and School of Mathematics, Institute for studies in Theoretical
Phisics and Mathematics (IPM), P.O. Box 19395-5746, Tehran, Iran.}
\email{naghipour@ipm.ir} \email {naghipour@tabrizu.ac.ir}
\email{d\_asadollahi@tabrizu.ac.ir}

\thanks{ 
This research was in part supported by a grant from IPM. \\
$^*$Corresponding author: e-mail: {\it naghipour@ipm.ir} (Reza
Naghipour)}%

\begin{abstract}
 Let $R$  denote a commutative Noetherian ring.  Brodmann et al. in \cite{BRS} defined and studied the concept of the local-global principle for annihilation
 of local cohomology modules at level $r\in\mathbb{N}$ for the ideals $\frak a$ and $\frak b$ of $R$.
It was shown in \cite{BRS} that this principle holds at levels 1,2, over $R$ and at all  levels whenever $\dim R\leq 4$. The goal of this paper is to show that,
if  the set $\Ass_R(H_{\fa}^{f_{\fa}^{\fb}(M)}(M))$ is finite or $f_{\fa}(M)\neq c_{\fa}^{\fb}(M)$, then the  local-global principle holds at all levels $r\in\mathbb{N}_0$, for all ideals $\fa , \fb$ of $R$ and each finitely generated $R$-module  $M$,  where $c_{\fa}^{\fb}(M)$ denotes the first non $\fb$-cofiniteness of local cohomology module
$H^i_{\fa}(M)$. As a consequence of this, we provide  a new and  short proof of the  Faltings' local-global principle for finiteness dimensions.
Also, several new results concerning the finiteness dimensions are given.
\end{abstract}
\maketitle
\section{Introduction}
Throughout this paper, let $R$ denote a commutative Noetherian ring
(with identity) and $\frak a$ an ideal of $R$. For an $R$-module $M$, the
$i^{\rm th}$ local cohomology module of $M$ with support in $V(\frak a)$
is defined as:
$$H^i_{\frak a}(M) = \underset{n\geq1} {\varinjlim}\,\, \Ext^i_R(R/\frak a^n, M).$$  We refer the reader to \cite{BS} or \cite{Gr1} for more
details about local cohomology.  An important theorem in local
cohomology is Faltings' Local-global Principle for the Finiteness
Dimension of local cohomology modules \cite[Satz 1]{Fa1}, which
states that for a positive integer $r$,  the $R_{\frak p}$-module $H^i_{\frak a R_{\frak p}}(M_{\frak p})$ is finitely generated
for all $i\leq r$ and for all ${\frak p}\in  \Spec R$ if and only if
the $R$-module $H^i_{\frak a}(M)$ is finitely generated for all $i\leq r$.

Another formulation of the Faltings' local-global principal,
particularly relevant for this paper, is in terms  of the finiteness dimension $f_{\frak a}(M)$ of $M$ relative
to $\frak a$, where
$$f_{\frak a}(M):=\inf\{i\in \Bbb{N}_0\,\,|\,\,H^i_{\frak a}(M)\,\,{\rm is}\,\,{\rm not}\,\,{\rm finitely}\,\,{\rm generated}
\}.\,\,\,\,\,\,\,\,\,\,\,\,\,\,\,\,\,\,\,\,\,\,\,\,\,(\dag)$$
 With the usual convention that the infimum of the empty set of integers is interpreted as $\infty$.
We can restate the Faltings' local-global principal $(\dag)$, in the following form
 $$f_{{\frak a}R_{\frak p}}(M_{\frak p})>r \,\,\,\,\, \text{ for all } {\frak p}\in \Spec(R) \Longleftrightarrow f_{\frak a}(M)>r. \,\,\,\,\,\,\,\,\,\,\,\,\,\,\,\,\,\,\,\,\,\,\,\,\,(\dag\dag)$$

Now, let $\frak b$ be a second ideal of $R$.  Recall that the $\frak b$-{\it finiteness dimension} of $M$ relative to $\frak a$ is defined by
\begin{eqnarray*}
f_{\frak a}^{\frak b}(M)&:=& \inf\{i\in \Bbb{N}_0\,\,|\,\,\frak b\not\subseteq {\rm Rad}(0:_{R}H^i_{\frak a}(M))\}\\
&=& \inf\{i\in \Bbb{N}_0\,\,|\,\,{\frak b}^nH^i_{\frak a}(M)\neq 0\,\,{\rm for}\,\,{\rm all}\,\,n\in \Bbb{N}\}.
\end{eqnarray*}

So it is rather natural to ask whether the Faltings' local-global principal, as stated in  $(\dag\dag)$, generalizes in the obvious way to the
invariants $f_{\frak a}^{\frak b}(M)$. In other words, is the statement
$$f_{{\frak a}R_{\frak p}}^{{\frak b}R_{\frak p}}(M_{\frak p})>r \,\,\,\,\, \text{ for all } {\frak p}\in{\Spec}(R) \Longleftrightarrow f_{\frak a}^{\frak b}(M)>r\,\,\,\,\,\,\,\,\,\,\,\,\,\,\,\,\,\,\,\,\,\,\,\,\,(\dag\dag\dag)$$
 true for each integer $r>0$? We shall say that  the local-global principle for the annihilation of local cohomology modules holds at level $r$ (over the ring $R$) if  $(\dag\dag\dag)$  is true (for the given $r$) for every choice of ideals ${\frak a}$, ${\frak b}$ of $R$ and every choice of finitely generated $R$-module $M$.

In fact, recently  Brodmann et al. in \cite{BRS} defined and studied the concept of the local-global principle for annihilation
 of local cohomology modules at level $r\in\mathbb{N}$ for the ideals $\frak a$ and $\frak b$ of $R$.
It is shown in \cite{BRS} that the local-global principle for the annihilation of local cohomology modules holds at levels 1,2, over an arbitrary commutative Noetherian ring $R$ and at all  levels whenever $\dim R\leq 4$.

More recently, for an non-negative integer $n$, Doustimehr and Naghipour in \cite{DN}, defined the $n$th $\fb$-finiteness dimension
of $M$ relative to $\fa$ (resp. the $n$th $\fb$-minimum $\fa$-adjusted depth of $M$) by
$$f_{\fa}^{\fb}(M)_n={\rm inf}\{i\in\mathbb{N}_0|\,\dim\Supp\fb^tH^i_{\fa}(M)\geq n\,\, \text{for all}\,\, t\in\mathbb{N}_0 \}.$$
(resp. $\lambda_{\fa}^{\fb}(M)_n={\rm inf}\{\lambda_{\fa R_{\p}}^{\fb R_{\p}}(M_{\p})|\,\,\p\in\Spec(R),\dim R/\p\geq n\}),$
and using the theory of G-dimension, gave a nice generalization of Faltings' Annihilator theorem.

Now, for a non-negative integer $n$, we define the upper $n$th $\fb$-finiteness dimension of $M$ relative to $\fa$ by
$$f_{\fa}^{\fb}(M)^n={\rm inf}\{f_{\fa R_{\p}}^{\fb R_{\p}}(M_{\p})|\,\, \p\in\Spec(R),\dim(R/\p)\geq n\}.$$

Note that $f_{\fa}^{\fb}(M)^n$ is either a positive integer or $\infty$, and
$$f_{\fa}^{\fa}(M)^n=f_{\fa}^n(M):=\inf\{f_{\fa R_{\frak p}}(M_{\frak p})\,\,|\,\,{\frak p}\in \Supp(M/\fa M)\,\,{\rm and}\,\,\dim R/{\frak p}\geq n\},$$
is the $n$th finiteness dimension of $M$ relative to $\fa$ (cf. \cite{BNS1}).

The purpose of the present paper is to show that, if the set $\Ass_R(H_{\fa}^{f_{\fa}^{\fb}(M)}(M))$ is finite or $f_{\fa}(M)\neq c_{\fa}^{\fb}(M)$, then the  Faltings' local-global principle holds at all levels $r\in\mathbb{N}_0$, where $c_{\fa}^{\fb}(M)$ denotes the first non $\fb$-cofiniteness of local cohomology module $H^i_{\fa}(M)$.  As a consequence of this, we provide a short proof of the  Faltings' local-global principle for finiteness dimensions. More precisely,
as a main result in the second section, we shall show that

\begin{thm}
Let $\fa$ and $\fb$ be ideals of $R$ such that $\fb\subseteq\fa$. Let  $M$ be a finitely generated $R$-module such that the set $\Ass_R(H_{\fa}^{f_{\fa}^{\fb}(M)}(M))$ is finite.  Then
$$f_{\fa}^{\fb}(M)={\rm inf}\{f_{\fa R_{\p}}^{\fb R_{\p}}(M_{\p})|\,\,\,\p\in\Spec(R)\}.$$
In particular,
$$f_{\fa}(M)={\rm inf}\{f_{\fa R_{\p}}(M_{\p})|\,\,\,\p\in\Spec(R)\}.$$
\end{thm}
The result in Theorem 1.1 is proved in Theorem 2.3 and Corollary 2.4. One of our tools for proving Theorem 1.1 is the following, which plays a key role in this paper.

\begin{prop}
Let $\fa$ and $\fb$ be ideals of $R$ such that $\fb\subseteq\fa$. Let  $M$ be a finitely generated $R$-module and $n$  a non-negative integer such that
 $(\Ass_R H_{\fa}^{f_{\fa}^{\fb}(M)_n}(M))_{\geq n}$ is finite. Then
$$f_{\fa}^{\fb}(M)_n=f_{\fa}^{\fb}(M)^n.$$
Here, for any subset $T$ of $\Spec(R)$ and a non-negative integer $n$, we set $$(T)_{\geq n}=\{ \p\in T | \dim R/\p\geq n\}.$$
\end{prop}

For a  finitely generated $R$-module $M$ and for ideals $\fa,\fb$  of $R$ with $\fb\subseteq\fa$,  the $\fb$-{\it cofiniteness dimension $c_{\fa}^{\fb}(M)$ of $M$ relative to} $\fa$ is define  by
$c_{\fa}^{\fb}(M)={\rm inf}\{i\in\mathbb{N}_0|\,H^i_{\fa}(M)\,\text{is not}\, \fb\text{-cofinite}\}.$
The notion of $\fb$-cofiniteness dimension $c_{\fa}^{\fb}(M)$ of $M$ relative to $\fa$ is introduced  and studied in \cite{BNS2}.

In Section 3, it is shown that  if   $f_{\fa}(M)\neq c_{\fa}^{\fb}(M)$, then the  local-global principle holds at all levels $r\in\mathbb{N}_0$, for all ideals $\fa , \fb$ of $R$ and each finitely generated $R$-module  $M$. Moreover,  we obtain some new results about the finiteness dimensions of local cohomology modules. In this section among other things, we derive
the following consequence of Theorem 1.1, which shows that the Faltings' local-global principle holds at all levels.

\begin{thm}
Let $\fb\subseteq\fa$ be ideals of $R$  and let  $M$ be a finitely generated $R$-module such that $f_{\fa}(M)\neq c_{\fa}^{\fb}(M)$.  Then the
Faltings' local-global principle for the finiteness of local cohomology modules  holds at all levels.
\end{thm}

 Throughout this paper, $R$ will always be a commutative Noetherian
ring with non-zero identity and $\frak a$ will be an ideal of $R$.  Recall
that an $R$-module $M$ is called $\frak a$-{\it cofinite} if $\Supp M\subseteq V(\frak a)$ and ${\rm Ext}^i_R(R/\frak a, M)$ is finitely
generated for all $i\geq 0$. The concept of $\frak a$-cofinite modules were introduced by Hartshorne \cite{Ha}. Also, if  $n$ is a non-negative integer,
then  $M$ is said to be in dimension $< n$, if there is a finitely generated submodule $N$ of $M$ such that $\dim \Supp  M/ N <n$ (cf. \cite[Definition 2.1]{AN}); and moreover if $T$ is a subset of $\Spec(R)$, then we define $$(T)_{\geq n}=\{ \p\in T | \dim R/\p\geq n\}.$$
 Finally, for any ideal $\frak{b}$ of $R$, the {\it radical} of $\frak{b}$, denoted by $\Rad(\frak{b})$, is defined to
be the set $\{x\in R \,: \, x^n \in \frak{b}$ for some $n \in \mathbb{N}\}$. For any unexplained notation and terminology we refer
the reader to \cite{BS} and \cite{Mat}.

\section{ Faltings' local-global principle and finiteness}

In this section we establish that the  Faltings' local-global principle for the finiteness of local cohomology modules holds at all levels over an arbitrary commutative Noetherian ring $R$,  whenever the set  $\Ass_R(H_{\fa}^{f_{\fa}^{\fb}(M)}(M))$ is finite. As a consequence,  we will provide  a new and  short proof of the  Faltings' local-global principle for finiteness dimensions. We begin with the following lemma which is needed in the proof of Proposition 2.2.
\begin{lem}
Let $\fa$ and $\fb$ be ideals of $R$ such that $\fb\subseteq\fa$, and let $M$ be a finitely generated $R$-module. Then, for every non-negative integer $n$,
$$f_{\fa}^{\fb}(M)_n\leq f_{\fa}^{\fb}(M)^n\leq \lambda_{\fa}^{\fb}(M)_n.$$
\end{lem}
\proof
Let $t:=f_{\fa}^{\fb}(M)^n$ and suppose that $t<f_{\fa}^{\fb}(M)_n$. Then, in view of the definition, there exists an integer $m$ such that $\dim\Supp (\fb^mH^{t}_{\fa}(M))<n$. Hence, for every prime ideal $\p$ of $R$ with $\dim R/\p\geq n$, we have
 $$(\fb^mH^{t}_{\fa}(M))_{\p}=0,$$
and so   $(\fb R_{\p})^mH^{t}_{\fa R_{\p}}(M_{\p})=0$, which is a contradiction.

Therefore $f_{\fa}^{\fb}(M)_n\leq f_{\fa}^{\fb}(M)^n$. The inequality $f_{\fa}^{\fb}(M)^n\leq \lambda_{\fa}^{\fb}(M)_n$, follows easily from the definition and  \cite[Theorem 9.3.5]{BS}. \qed\\

The following proposition  plays a key role in the proof of the main result of this section.
\begin{prop}
Let $\fa$ and $\fb$ be ideals of $R$ such that $\fb\subseteq\fa$. Let  $M$ be a finitely generated $R$-module and  let $n$ be a non-negative integer such that
 $(\Ass_R H_{\fa}^{f_{\fa}^{\fb}(M)_n}(M))_{\geq n}$ is finite. Then
$$f_{\fa}^{\fb}(M)_n=f_{\fa}^{\fb}(M)^n.$$
\end{prop}
\proof
Let $t:=f_{\fa}^{\fb}(M)_n$ and suppose that
$$(\Ass_R H_{\fa}^t(M))_{\geq n}=\{\p_1,\dots,\p_r\}.$$
In view of Lemma 2.1, it is enough for us to show that $f_{\fa}^{\fb}(M)^n\leqslant t$. Suppose the contrary that $t<f_{\fa}^{\fb}(M)^n$, and look for a contradiction.
To this end, in view of definition  $t<f_{\fa R_{\p_i}}^{\fb R_{\p_i}}(M_{\p_i})$, for all $1\leq i\leq r$. Hence, for all $1\leq i\leq r$, there exists an integer $n_i$ such that
$$(\fb R_{\p_i})^{n_i}H^{t}_{\fa R_{\p_i}}(M_{\p_i})=0.$$
Set $n:={\rm max}\{n_1,\ldots,n_r\}$.  Then, for all $1\leq i\leq r$, we have
$$(\fb^nH^{t}_{\fa}(M))_{\p_i}=0.$$
Now, it is easy to see that $\dim\Supp(\fb^nH^{t}_{\fa}(M))<n$, which is a contradiction.\qed\\

We are now ready to state and prove the main theorem of this section, which shows
that the Faltings' local-global principle for the finiteness of local cohomology modules is
valid at all levels over any Noetherian ring $R$,  whenever the set  $\Ass_R(H_{\fa}^{f_{\fa}^{\fb}(M)}(M))$ is finite.
\begin{thm}
Let $\fa$ and $\fb$ be ideals of $R$ such that $\fb\subseteq\fa$. Let  $M$ be a finitely generated $R$-module such that the set $\Ass_R(H_{\fa}^{f_{\fa}^{\fb}(M)}(M))$ is finite. Then
$$f_{\fa}^{\fb}(M)={\rm inf}\{f_{\fa R_{\p}}^{\fb R_{\p}}(M_{\p})|\,\,\,\p\in\Spec(R)\}.$$
\end{thm}
\proof
Put $n=0$ in the Proposition 2.2, and use the fact that $f_{\fa}^{\fb}(M)=f_{\fa}^{\fb}(M)_0$.\qed \\

An immediate  consequence of Theorem 2.3, we give a short proof of the  Faltings' local-global principle for finiteness dimensions.

\begin{cor}{\rm (\cite[\bf 9.6.2 Local-global Principle for Finiteness Dimensions]{BS})}
Let $M$ be a finitely generated $R$-module. Then for  any ideal $\fa$ of $R$,
$$f_{\fa}(M)={\rm inf}\{f_{\fa R_{\p}}(M_{\p})|\,\,\,\p\in\Spec(R)\}.$$
\end{cor}
\proof
Put $\fa=\fb$ in Theorem 2.3 and use \cite[Corollary 2.3 ]{BL}.\qed\\

\section{Faltings' local-global principle and cofiniteness}

The purpose of this section  is to show that  if  $f_{\fa}(M)\neq c_{\fa}^{\fb}(M)$, then the  local-global principle holds at all levels $r\in\mathbb{N}_0$, for all ideals $\fa , \fb$ of $R$ and each finitely generated $R$-module  $M$. Moreover,  we obtain some new results about the finiteness dimensions of local cohomology modules. In order to do this, let us recall that for a  finitely generated $R$-module $M$ and for ideals $\fa,\fb$  of $R$ such that $\fb\subseteq\fa$,  the $\fb$-cofiniteness dimension $c_{\fa}^{\fb}(M)$ of $M$ relative to $\fa$ is defined  by
$c_{\fa}^{\fb}(M)={\rm inf}\{i\in\mathbb{N}_0|\,H^i_{\fa}(M)\,\text{is not}\, \fb\text{-cofinite}\}.$
The notion of $\fb$-cofiniteness dimension $c_{\fa}^{\fb}(M)$ of $M$ relative to $\fa$ is introduced  and studied in \cite{BNS2}.
\begin{prop}
Let $\fa$ and $\fb$ be ideals of $R$ such that $\fb\subseteq\fa$  and let  $M$ be a finitely generated $R$-module. Then
$$f_{\fa}(M)={\rm min}\{f_{\fa}^{\fb}(M),c_{\fa}^{\fb}(M)\}.$$
\end{prop}
\proof
It is clear that $f_{\fa}(M)\leq f_{\fa}^{\fb}(M)$. Now, if $t:=c_{\fa}^{\fb}(M)<f_{\fa}(M)$, then $H^t_{\fa}(M)$ is a finitely generated $R$-module. Since $$\Supp H^t_{\fa}(M)\subseteq V(\fa)\subseteq V(\fb),$$
it follows that the $R$-module $H^t_{\fa}(M)$ is $\fb$-cofinite, which is a contradiction. Whence $f_{\fa}(M)\leq c_{\fa}^{\fb}(M)$ and so
$$f_{\fa}(M)\leq{\rm min}\{f_{\fa}^{\fb}(M),c_{\fa}^{\fb}(M)\}.$$
Finally, suppose that
$$r:=f_{\fa}(M)<{\rm min}\{f_{\fa}^{\fb}(M),c_{\fa}^{\fb}(M)\},$$ and look for a contradiction.
To do this, as $r<f_{\fa}^{\fb}(M)$, there exists an integer $n$ such that $\fb^nH^r_{\fa}(M)=0$, and thus
$$H^r_{\fa}(M)\cong {\rm Hom}_R(R/\fb^n,H^r_{\fa}(M)).$$
Now, since $r<c_{\fa}^{\fb}(M)$,  it follows that the $R$-module $H^r_{\fa}(M)$ is $\fb$-cofinite, and so it yields from the above isomorphism that
$H^r_{\fa}(M)$ is finitely generated, which is a contradiction. \qed\\

We are now ready to state and prove the main theorem of this section, which shows
that the Faltings' local-global principle for the finiteness of local cohomology modules is
valid at all levels over any Noetherian ring $R$,  whenever $f_{\fa}(M)\neq c_{\fa}^{\fb}(M)$.

\begin{thm}
Let $\fb\subseteq\fa$ be ideals of $R$  and let  $M$ be a finitely generated $R$-module such that $f_{\fa}(M)\neq c_{\fa}^{\fb}(M)$.  Then the
Faltings' local-global principle for the finiteness of local cohomology modules  holds at all levels.
\end{thm}
\proof
Since $f_{\fa}(M)\neq c_{\fa}^{\fb}(M)$, it follows from Proposition 3.1 that $f_{\fa}^{\fb}(M)\leq c_{\fa}^{\fb}(M)$, and so
$f_{\fa}(M)= f_{\fa}^{\fb}(M)$.  Now, the assertion follows from Theorem 2.3 and Brodmann-Lashgari's result \cite[Corollary 2.3 ]{BL}.  \qed\\

Before bringing the next result, we give a couple of lemmas that will be used in the proof of Theorem 3.5.

\begin{lem}
Let  $M$ be an  $R$-module and let $s$  be a non-negative integer. Then

 ${\rm(i)}$ If $M$ is in dimension $<s$, then $M_{\frak p}$ is a finitely generated $R_{\frak p}$-module for all $\frak p \in \Spec (R)$ with $\dim R/\frak p\geq s$.

${\rm(ii)}$  If $M$ is in dimension $<s$, then the set $(\Ass_R M)_{\geq s}$ is finite.
\end{lem}
\proof
The part $\rm(i)$ follows easily from the definition. In order to prove $\rm(ii)$,  since $M$ is in dimension $< s$, it follows from the definition that there is a finitely generated submodule $M'$ of $M$ such that $\dim \Supp M/M'<s$. Now, from the exact sequence $$0 \longrightarrow M' \longrightarrow M\longrightarrow M/M' \longrightarrow 0$$
we obtain
$$(\Ass_R M)_{\geq s}\subseteq (\Ass_RM')_{\geq s}\cup (\Ass_R M/M')_{\geq s}.$$
As $\dim \Supp M/M'<s$, it follows that $(\Ass_R M/M')_{\geq s}=\emptyset$,  and so  the set $(\Ass_R M)_{\geq s}$ is finite.\qed \\

Before bringing the next lemma, let us recall that a full subcategory $\mathcal{S}$ of the category of
$R$-modules is called a {\it Serre subcategory}, when it is closed under taking submodules, quotients
and extensions.
\begin{lem}
For any non-negative integer $n$, the class of in dimension $<n$ modules over a  Noetherian ring $R$ consists  a  Serre subcategory of
the category of $R$-modules.
\end{lem}

\proof See \cite[Corollary 2.3]{MNS}.  \qed \\

\begin{thm}
Let $(R,{\frak m})$ be a complete local ring and $M$ a finitely generated $R$-module. Let  $\fb\subseteq\fa$ be ideals of $R$  and let $n$  be a non-negative
integer such  that $f_{\fa}^{\fb}(M)_n=f_{\fa}^{\fb}(M)^n$.
Then $f_{\fa}^{\fb}(M)_n\geq {\rm depth}(\fb,M)$ if and only if $f_{\fa}^n(M)\geq {\rm depth}(\fb,M)$.
\end{thm}

\proof
Let ${\rm depth}(\fb,M)=s$.  Since $\fb\subseteq\fa$, it yields that  $f_{\fa}(M)\leq f_{\fa}^{\fb}(M)$ and so
$$f_{\fa}^n(M)\leq f_{\fa}^{\fb}(M)^n=f_{\fa}^{\fb}(M)_n.$$
Recall that $f_{\fa}^n(M):=\inf\{f_{\fa R_{\frak p}}(M_{\frak p})\,\,|\,\,{\frak p}\in \Supp(M/\fa M)\,\,{\rm and}\,\,\dim R/{\frak p}\geq n\}.$

Now, suppose that $f_{\fa}^{\fb}(M)_n\geq s$ and we show that $f_{\fa}^n(M)\geq s$.  We use induction on $s$. For $s=0$ there is nothing to show. So assume that $s>0$ and the result has been proved for $s-1$.

Since ${\rm depth}(\fb,M)>0$, it follows that $\fb$ contains an element $x$ which is a non-zero divisor on $M$.
Moreover, as $f_{\fa}^{\fb}(M)_n\geq s$, there exists an integer $t$ such that, for all $i<s$,  $\dim\Supp(\fb^t H^{i}_{\fa}(M))<n$.
Hence from the exact sequence
$$0\longrightarrow M \stackrel{x^t}\longrightarrow M\longrightarrow M/x^tM\longrightarrow 0,$$
and \cite[Lemma  2.9]{DN}, there exists an integer $l$ such that, for all $i<s-1$, $$\dim\Supp(\fb^lH^{i}_{\fa}(M/x^tM))<n,$$  and thus $f_{\fa}^{\fb}(M/x^tM)_n\geq s-1$.
Therefore, in view of the inductive hypothesis $$f_{\fa}^n(M/x^tM)\geq s-1,$$  and so by virtue of \cite[Theorem 2.5]{AN}, $H^{s-2}_{\fa}(M/x^tM)$ is in dimension $<n$.
Therefore, it follows from the exact sequence and Lemma 3.4
$$0\longrightarrow H^{s-2}_{\fa}(M)/x^tH^{s-2}_{\fa}(M) \longrightarrow H^{s-2}_{\fa}(M/x^tM)\longrightarrow(0:_{H^{s-1}_{\fa}(M)}x^t)\longrightarrow 0,$$
that $(0:_{H^{s-1}_{\fa}(M)}x^t)$ is in dimension $<n$. Hence in view of Lemma 3.3, the $R_{\frak p}$-module
$(0:_{(H^{s-1}_{\fa}(M))_{\p}}x^t/1)$ is finitely generated, for all $\frak p \in \Spec(R)$ with  $\dim R/\p\geq n$.  Now, since $\dim\Supp(\fb^t H^{s-1}_{\fa}(M))<n$, we have $(\fb^t H^{s-1}_{\fa}(M))_{\p}=0$, and hence
 $x^t/1 (H^{s-1}_{\fa}(M))_{\p}=0.$
 Consequently,
 $$(H^{s-1}_{\fa}(M))_{\p}=(0:_{(H^{s-1}_{\fa}(M))_{\p}}x^t/1),$$
is a finitely generated $R_{\frak p}$-module. Therefore, by virtue of Lemma 3.3,   the $R$-module  $ H^{s-1}_{\fa}(M)$ is also in dimension $<n$.
 Consequently,  we deduce from \cite[Theorem 2.5]{AN} that $f_{\fa}^n(M)\geq s$, as required.\qed\\

 \begin{cor}
Let $(R,{\frak m})$ be a complete local ring and $M$ a finitely generated $R$-module. Let  $\fb\subseteq\fa$ be ideals of $R$  and let $n$  be a non-negative
integer such  that ${\rm depth}(\fb,M)\leq f_{\fa}^{\fb}(M)_n$ and  the set $\Ass_R(H_{\fa}^{f_{\fa}^{\fb}(M)}(M))$ is finite.
Then $${\rm depth}(\fb,M)\leq \min\{f_{\fa}^{\fb}(M)_n,  f_{\fa}^n(M)\}.$$
\end{cor}
\proof
The result follows from Proposition 2.2 and Theorem 3.5. \qed\\

 \begin{thm}
Let $(R,{\frak m})$ be a complete local ring and $M$ a finitely generated $R$-module. Let  $\fa$ be an ideal of $R$  and  let $s,t$ be   non-negative integers.
Then the following statements are equivalent.

$\rm(i)$  $H^{i}_{\fa}(M)$ is in dimension $<s$, for all $i<t$.

$\rm(ii)$ There exists an integer $n$ such that $\dim\Supp(\fa^nH^{i}_{\fa}(M))<s$, for all $i<t$.
\end{thm}
\proof
In order to show $\rm(i)\Longrightarrow \rm(ii)$, let $i$ be an integer such that $i<t$.  Then, since $H^{i}_{\fa}(M)$ is in dimension $<s$, in view of Lemma 3.3, the set
$$(\Ass_R H^{i}_{\fa}(M))_{\geq s}$$
is finite. Let
$$(\Ass_R H^{i}_{\fa}(M))_{\geq s}=\{\p_1,\dots,\p_n\}.$$
Thus, for all $j$ with $1\leq j\leq r$, the $R_{\p_j}$-module $(H^{i}_{\fa}(M))_{\p_j}$ is  finitely generated
and thus there exists $n_j\in\mathbb{N}$ such that $(\fa^{n_j}H^{i}_{\fa}(M))_{\p_j}=0$.

Put $n:=\max\{n_1,\ldots,n_r\}$. Then,  for all $j$ with $1\leq j\leq r$, we have $(\fa^{n}H^{i}_{\fa}(M))_{\p_j}=0$.
Now, we show that $\dim\Supp(\fa^{n}H^{i}_{\fa}(M))<s$.
For this, let $\q\in\Ass_R(\fa^{n}H^{i}_{\fa}(M))$.  Then, if $\dim R/\q\geq s$, then there exists $j$ such that $\q=\p_j$ and this is a contradiction. Thus $\dim R/\q< s$,  as desired.

Now, we show $\rm(ii)\Longrightarrow \rm(i)$.  To do this, for all $i<t$ and for all prime ideals $\p$  with $\dim R/\p\geq s$,
we have $(\fa^{n}H^{i}_{\fa}(M))_{\p}=0$.
 Thus, for all $i<t$ and for all prime ideals $\p$ with $\dim R/\p\geq s$, the $R_{\p}$-module $ H^{i}_{\fa}(M))_{\p}$ is  finitely generated. Hence $f_{\fa}^s(M)\geq t$. Now, since $$f_{\fa}^s(M)={\rm inf} \{0\leq i\in\inte|\,H^{i}_{\fa}(M)\text{ is not in dimension}<s\},$$
(cf.  \cite[Theorem 2.5]{AN}), it follows that $H^{i}_{\fa}(M)$ is in dimension $<s$ for all $i<t$, as required.\qed\\

\begin{cor}
Let $(R,{\frak m})$ be a complete local ring,  $M$ a finitely generated $R$-module and $\fa$  an ideal of $R$. Then
for any non-negative integer $n$,  $$f_{\fa}^n(M)=f_{\fa}^{\fa}(M)_n.$$
\end{cor}

 \proof
 The assertion follows from Theorem 3.7 and  the definitions of  $f_{\fa}^n(M)$ and $f_{\fa}^{\fa}(M)_n.$  \qed\\

\begin{thm}
Let $R$ be a Noetherian ring,  $M$ a finitely generated $R$-module and let  $\fb\subseteq\fa$ be  ideals of $R$. Then, for any non-negative integers
$i$ and $ s$,  the following statements are equivalent.

$\rm (i)$  There exists an integer $n$ such that $\dim\Supp(\fb^nH^{i}_{\fa}(M))<s$.

$\rm (ii)$  There exists an integer $m$ such that $\fb^mH^{i}_{\fa}(M)$ is in dimension $<s$.
\end{thm}

\proof
The implication $\rm (i) \Longrightarrow \rm (ii)$ is clear. In order to show  $\rm (ii) \Longrightarrow \rm (i)$,
since $\fb^mH^{i}_{\fa}(M)$ is in dimension $<s$,  it follows from Lemma 3.3 that the set
$$(\Ass_R \fb^mH^{i}_{\fa}(M))_{\geq s}$$
is finite.  Let
$$(\Ass_R \fb^mH^{i}_{\fa}(M))_{\geq s}=\{\p_1,\dots,\p_r\}.$$
Then, for all $j$ with $1\leq j\leq r$, the $R_{\p_j}$-module $(\fb^mH^{i}_{\fa}(M))_{\p_j}$ is finitely generated. Hence,  there exists $n_j\in\mathbb{N}$ such that $(\fa^{n_j}(\fb^mH^{i}_{\fa}(M))_{\p_j}=0$. Since $\fb\subseteq\fa$, it follows that $(\fb^{m+n_j}H^{i}_{\fa}(M))_{\p_j}=0$.

Set $n:=\max\{m+n_1,\ldots,m+n_r\}$. Then,  for all $j$ with $1\leq j\leq r$, we deduce that $(\fb^{n}H^{i}_{\fa}(M))_{\p_j}=0$.

We  show that $\dim\Supp(\fb^{n}H^{i}_{\fa}(M))<s$.  To do this, let $\p\in\Ass_R(\fb^{n}H^{i}_{\fa}(M))$. Then $\p\in \Ass_R(\fb^{m}H^{i}_{\fa}(M))$.
Now, if $\dim R/\p\geq s$, then there exists $j$ such that $\p=\p_j$ and this is a contradiction. Thus $\dim R/\p< s$,  as required.\qed\\

\begin{cor}
Let $R$ be a Noetherian ring and $M$ a finitely generated $R$-module.
Let $\fa$ and $\fb$ be two  ideals of $R$ such that $\fb\subseteq\fa$. Then, for any non-negative integer
$n$,
$$f_{\fa}^{\fb}(M)_n={\rm inf} \{0\leq i\in\inte|\,\fb^m H^{i}_{\fa}(M)\text{ is not in dimension} <n\,\, \text{for all}\,\, m\}.$$
\end{cor}

\proof
The assertion follows from Theorem 3.9 and the definition of $f_{\fa}^{\fb}(M)_n$. \qed\\

\begin{center}
{\bf Acknowledgments}
\end{center}
The authors would like to thank Professors  Hossein Zakeri and Kamal Bahmanpour for reading of the first draft and valuable discussions.
Also, we  would like to thank from the Institute for Research in Fundamental Sciences (IPM) for the financial support.

\end{document}